\documentclass[a4paper, 12pt]{article}
\usepackage{enumerate, theorem}
\usepackage{amsmath, amsfonts, amssymb}
\usepackage[height=22.5cm, width=15cm]{geometry}
\usepackage[all]{xy}
\usepackage{mathrsfs, yfonts}

\DeclareMathOperator{\id}{id}

\DeclareMathOperator{\C}{\mathbb{C}}

\newcommand{\A}{\tilde{\mathcal{A}}}

\newcommand{\parag}[1]{\paragraph{\sc{#1.}}}

\newtheorem{thm}{Th\'eor\`eme}[subsection]
\newtheorem{defn}[thm]{D\'efinition}
\newtheorem{cor}[thm]{Corollaire}
\newtheorem{prop}[thm]{Proposition}
\newtheorem{lemma}[thm]{Lemme}

\begin{document}

\title{Changements de variable pour un th\`eme.}

\author{Daniel Barlet}

\date{17/03/10}
\maketitle

\section*{Abstract}
We study the behaviour of the notion of "thema", introduced in our previous article [B.09b], by a change of variable. We show not only that the fundamental invariants of such a thema, corresponding to the Bernstein polynomial, are stable by a change of variable, but also other numerical invariants called principal parameters. \\
We show on a rank 3  example that nevertheless the isomorphism class of a thema is not stable in general by a change of variable. We conclude in proving that a change of variable transforms an holomorphic family of thema in an holomorphic family. This implies that  non principal parameters change holomorphically.

\tableofcontents

\section{Introduction}

Quand on consid\`ere un morphisme propre \ $\tilde{f} : X \to C$ \ d'une vari\'et\'e complexe \ $X$ \ sur une courbe lisse \ $C$, telle que l'on ait \ $\{ d\tilde{f } = 0 \} \subset \tilde{f}^{-1}(s_0)$ \ au voisinage du point \ $s_0$ \ de \  \ $ C$, le choix d'une coordonn\'ee locale sur \ $C$ \ pr\`es de \ $s_0$ \ permet de se ramener \`a la situation de d\'eg\'enescence standard d'une famille de vari\'et\'es complexes param\`etr\'ee par un disque  \ $f : X \to D \subset \C$, la fibre singuli\`ere \'etant au-dessus de l'origine. On est alors dans la situation o\`u l'on dispose d'une fonction holomorphe propre \ $f$ \  et on peut alors construire les faisceaux de (a,b)-modules d\'eduits des complexes \ $(Ker\, df)^{\bullet}, d^{\bullet})$ \ gr\^ace \`a  la multiplication par la fonction num\'erique \ $f$ \ (voir [B.08] th\'eor\`eme 2.1.1).\\
Mais il est clair que la construction, et  pr\'ecis\'ement  via la fonction num\'erique \ $f$,  d\'epend du choix de la coordonn\'ee locale choisie pr\`es du point \ $s_0$ \ de la courbe \ $C$. L'objet de ce qui suit est d'\'etudier le comportement de nos constructions par un changement de coordonn\'ee locale centr\'ee en \ $s_0$. Le r\'esultat que nous obtenons est le th\'eor\`eme suivant.

\begin{thm}\label{chgt de variable}
Soit \ $E$ \ un th\`eme \ $[\lambda]-$primitif de rang \ $k$ \ et d'invariants fondamentaux \ $\lambda_1, p_1, \dots, p_{k-1}$. Soit \ $\theta(a) : = a + \theta_2.a^2 + \dots \in \C[[a]]$ \ un changement de variable. Alors le (a,b)-module \ $\theta_*(E)$ \ est un th\`eme  \ $[\lambda]-$primitif de rang \ $k$ \ ayant les m\^emes invariants fondamentaux que \ $E$. De plus les param\`etres principaux de \ $\theta_*(E)$ \ sont les m\^emes que ceux de \ $E$.
\end{thm}

Le lecteur trouvera les d\'efinitions des invariants fondamentaux et des param\`etres principaux d'un th\`eme \ $[\lambda]-$primitif dans les rappels de la section 2. Mais on notera d\'ej\`a que les invariants fondamentaux d'un th\`eme \ $[\lambda]-$primitif sont localement constants dans une famille holomorphe et d\'eterminent le polyn\^ome de Bernstein de ce th\`eme, alors que les param\`etres principaux sont des nombres complexes (non nuls) qui ne d\'ependent que de la classe d'isomorphisme du th\`eme consid\'er\'e, mais varient holomorphiquement avec le param\`etre dans une famille holomorphe. \\
Nous montrerons  par un exemple explicite, que, cependant, la classe d'isomorphisme d'un  th\`eme \ $[\lambda]-$primitif de rang \ $k \geq 3$ \ n'est pas, en g\'en\'eral, invariante par changement de variable, ce qui montre l'int\'er\^et de l'invariance des param\`etres fondamentaux donn\'ee par notre r\'esultat.\\
Nous conclueront en montrant que la notion de famille holomorphe de th\`emes \ $[\lambda]-$primitifs est stable par changement de variable. 

\bigskip

\section{Position du probl\`eme.}

\parag{Les alg\`ebres \ $\A$ \ et \ $\hat{A}$} Notons \ $\A$ \ la \ $\C-$alg\`ebre unitaire suivante :
$$ \A : = \{ \sum_{\nu = 0}^{\infty} \ P_{\nu}(a).b^{\nu} \} $$
o\`u les \ $P_{\nu}$ \ sont des polyn\^omes de \ $\C[x]$. La multiplication est d\'efinie par la relation de commutation \ $a.b - b.a = b^2$ \ et le fait que les multiplications \`a gauche et \`a droite par \ $a$ \ sont continues pour la topologie \ $b-$adique. On a alors pour chaque \ $S \in \C[[b]]$ \ la relation de commutation \ $a.S = S.a + b^2.S' $ \ o\`u \ $S'$ \ d\'esigne la d\'eriv\'ee de \ $S$ \ par rapport \`a la variable \ $b$. On v\'erifie facilement que l'on a \ $b^n.\A = \A.b^n$ \ et que \ $\A$ \ est compl\`ete pour la topologie \ $b-$adique.\\
Il est aussi facile de voir que l'\'el\'ement \ $1 + a \in \A$ \ n'est pas inversible dans \ $\A$. On d\'efinit l'alg\`ebre \ $\hat{A}$ \ comme la compl\'et\'ee \ $a-$adique de  \ $\A$ \ comme \ $\C[a]-$module \`a gauche, c'est-\`a-dire que
$$ \hat{A} : =  \{ \sum_{\nu = 0}^{\infty} \ P_{\nu}(a).b^{\nu} \} $$
o\`u maintenant les \ $P_{\nu}$ \ sont des s\'eries formelles dans \ $\C[[x]]$. Les relations de commutations dans \ $\A$ \ 
$$ a^n.b = b.a^n + n.b.a^{n-1}.b \quad \forall n \in \mathbb{N}$$ 
donnent alors, pour \ $T \in \C[[x]]$ \ la relation \ $T(a).b = b.T(a) + b.T'(a).b $ \ o\`u \ $T' \in \C[[x]]$ \ est la d\'eriv\'ee "usuelle" de la s\'erie formelle \ $T$.

\begin{defn}\label{a,b}
Un {\bf (a,b)-module} \ $E$ \ est un \ $\A-$module \`a gauche qui est libre de type fini sur  la sous-alg\`ebre \ $\C[[b]]$ \ de \ $\A$. Il sera {\bf \`a p\^ole simple} s'il v\'erifie \ $a.E \subset b.$. Il est {\bf r\'egulier}  s'il est contenu dans un (a,b)-module \`a p\^ole simple. Il est {\bf local} s'il existe \ $N \in \mathbb{N}$ \ tel que l'on ait \ $a^N.E \subset b.E$.
\end{defn}

Un (a,b)-module r\'egulier est toujours local mais la r\'eciproque est inexacte. Un (a,b)-module est, par d\'efinition, complet pour la topologie \ $b-$adique. S'il est local, il est \'egalement complet pour la topologie \ $a-$adique. C'est donc naturellement un \ $\hat{A}-$module \`a gauche.\\
En fait tous les (a,b)-modules qui vont nous int\'eresser \'etant r\'eguliers, ils sont naturellement des \ $\hat{A}-$modules \`a gauche.

Nous renvoyons le lecteur \`a [B.09a]  pour les d\'efinitions de base sur les (a,b)-modules r\'eguliers que nous utiliserons dans ce qui suit.

\parag{Exemple} Soit \ $\lambda \in \C, \Re(\lambda) > 0 $ \ et soit \ $k \geq 1$ \ un entier. Notons 
$$ \Xi_{\lambda}^{(k-1)} : = \sum_{j=0}^{k-1}  \C[[b]].s^{\lambda-1}.\frac{(Log\,s)^j}{j!} $$
le (a,b)-module d\'efini par les relations
\begin{align*}
& a.s^{\lambda-1} : = \lambda.b.s^{\lambda-1} \quad {\rm et} \\ 
& a.s^{\lambda-1}.\frac{(Log\,s)^j}{j!} = \lambda.b.s^{\lambda-1}.\frac{(Log\,s)^j}{j!} + b.\big[s^{\lambda-1}.\frac{(Log\,s)^{j-1}}{(j-1)!}\big] \quad {\rm pour} \quad j \in[1,k-1].
\end{align*}

 C'est clairement un (a,b)-module \`a p\^ole simple de rang \ $k$ \ sur \ $\C[[b]]$.\hfill $\square$\\
 
 \begin{defn}\label{theme primitif}
 Nous appellerons {\bf th\`eme \ $[\lambda]-$primitif} un sous-$\A-$module (\`a gauche) {\bf monog\`ene} d'un \ $ \Xi_{\lambda}^{(N-1)}$ \ pour un \ $N \in \mathbb{N}^*$ \ et un\ $\lambda \in \mathbb{Q}^{+*}$.
 \end{defn}
 
 Pour tout \ $\varphi \in  \Xi_{\lambda}^{(N-1)}$, le sous$-\A-$module \ $E : = \A.\varphi$ \ est un th\`eme. En effet c'est un sous-$\C[[b]]-$module d'un \ $\C[[b]]-$module libre de rang fini, il est donc libre de rang fini. Comme \ $ \Xi_{\lambda}^{(N-1)}$ \ est \`a p\^ole simple \ $E$ \ est r\'egulier. Et r\'eciproquement tout  th\`eme \ $[\lambda]-$primitif est de cette forme. \\
 Il est facile de voir que si \ $k-1$ \ est la puissance maximale de \ $Log\, s$ \ qui appara\^it effectivement dans \ $\varphi$ \ le rang de \ $E$ \ sur \ $\C[[b]]$ \ est \'egal \`a \ $k$, et l'on peut alors prendre \ $N = k$.
 
 \parag{Rappels} Il a \'et\'e d\'emontr\'e dans [B.09b] qu'un th\`eme \ $[\lambda]-$primitif \ $E$ \ admet une unique suite de Jordan-H{\"o}lder et que cette propri\'et\'e caract\'erise les th\`emes \ $[\lambda]-$primitifs parmi les (a,b)-modules monog\`enes g\'eom\'etriques\footnote{Un (a,b)-module r\'egulier est g\'eom\'etrique si les racines de son polyn\^ome de Bernstein sont dans \ $-\mathbb{Q}^{*+}$.}. Dans ce cas, si \ $k  = rg(E)$ \  on a pour chaque \ $j \in[0,k]$ \ un unique sous-(a,b)-module \ $F_j$ \ qui est normal de rang \ $j$.\\

\begin{defn}[Invariants fondamentaux]
La suite \ $\lambda_1, \dots, \lambda_k$ \  des exposants des quotients successifs de la suite de Jordan-H{\"o}lder  d'un th\`eme \ $[\lambda]-$primitif est \'equivalente \`a la donn\'ee de \ $\lambda_1, p_1, \dots, p_{k-1}$ \ o\`u les \ $p_j$ \ sont d\'efinis en posant \ $\lambda_{j+1} = \lambda_j + p_j - 1$. Les \ $p_j$ \ sont dans \ $\mathbb{N}$. Nous appellerons les nombres \ $\lambda_1, p_1, \dots, p_{k-1}$ \ {\bf les invariants fondamentaux} du th\`eme \ $[\lambda]-$primitif consid\'er\'e
\end{defn}

On notera que si \ $E \subset \Xi^{k-1}_{\lambda}$ \ avec \ $ k  = rg(E)$, ce qui est toujours possible d'apr\`es notre d\'efinition, on a \ $F_j = E \cap \Xi^{j-1}_{\lambda}$.\\

La classification des th\`emes de rang 1  est imm\'ediate, puisqu'elle se r\'eduit \`a la donn\'ee d'un rationnel positif \ $\lambda = \lambda_1 \in \mathbb{Q}^{+*}$. \\
Rappelons la classification des th\`emes \  $[\lambda]-$primitifs de rang 2 est donn\'ee par le th\'eor\`eme suivant qui se d\'eduit facilement de la classification des (a,b)-modules r\'eguliers de rang 2 donn\'ee dans [B.93] (voir [B.09b]).

\begin{thm}\label{rang 2}
Fixons \ $\lambda_1 > 1$ \ dans \ $\mathbb{Q}^{+*}$ \ et \ $p_1 \in \mathbb{N}$. Les classes d'isomorphismes de th\`emes \  $[\lambda]-$primitifs de rang 2 d'invariants fondamentaux est la suivante :
\begin{enumerate}
\item Pour \ $p_1 = 0$ \ on a un unique th\`eme \  $[\lambda]-$primitifs de rang 2 d'invariants fondamentaux \ $(\lambda_1, 0)$ \ \`a isomorphisme pr\`es ; il est donn\'e par le \ $\A-$module
$$ \A\big/\A.(a - \lambda_1.b)(a - (\lambda_1-1).b).$$
\item Pour chaque  \ $p \geq 1$ \ les classes d'isomorphismes de th\`emes  \  $[\lambda]-$primitifs de rang 2 d'invariants fondamentaux \ $(\lambda_1, p_1)$ \ sont en bijection avec \ $\C^*$. Au nombre \ $\alpha \in \C^*$ \ correspond la classe d'isomorphisme du \ $\A-$module (\`a gauche)
$$\A\big/\A.(a - \lambda_1.b)(1 + \alpha.b^{p_1})^{-1}(a - (\lambda+p_1-1).b) .$$
\end{enumerate}
\end{thm}

\begin{defn}\label{parametre}
Pour \ $p \geq 1$ \ le nombre \ $\alpha \in \C^*$ \ caract\'erisant la classe d'isomorphisme d'un th\`eme \ $E$ \ qui est   \  $[\lambda]-$primitif de rang 2 d'invariants fondamentaux \ $(\lambda_1, p_1)$ \ sera appel\'e le {\bf param\`etre} de \ $E$.\\
Dans le cas \ $p_1 = 0$ \ nous conviendont que le param\`etre est \ $\alpha : = \emptyset$.
\end{defn}

\begin{defn} Soit \ $E$ \ un th\`eme \ $[\lambda]-$primitif et soit \ $(F_j)_{j\in[0,k]}$ \ ses sous-th\`emes normaux param\'etr\'es par leurs rangs respectifs. Nous appellerons {\bf param\`etres principaux} du th\`eme \ $E$ \ la liste (ordonn\'ee) des param\`etres des th\`emes de rang 2 \ $F_{j+1}\big/F_{j-1}$ \ pour \ $j \in[1,k-1]$. On les notera \ $\alpha_1, \dots, \alpha_{k-1}$.
\end{defn}

 Ce sont des nombres qui ne d\'ependent que de la classe d'isomorphisme de \ $E$. Mais ils ne suffisent en g\'en\'eral pas \`a determiner la classe d'isomorphisme de \ $E$ \ en rang \ $\geq 3$. On notera que ces nombres sont toujours non nuls, et que,si \ $p_j \not= 0$ \  \ $\alpha_j$ \ peut prendre n'importe quelle valeur dans \ $\C^*$.

\bigskip

Il est d\'emontr\'e dans [B.09b] que les param\`etres  sont les coefficients \ $\alpha_1, \dots, \alpha_{k-1}$ \ respectivement de \ $b^{p_j}$ \ dans \ $S_j$ \ si l'annulateur d'un  g\'en\'erateur du th\`eme est l'id\'eal \`a gauche de \ $\A$ \ engendr\'e par l'\'el\'ement
$$ (a - \lambda_1.b).S_1^{-1}.(a- \lambda_2.b).S_2^{-1} \dots S_{k-1}^{-1}.(a - \lambda_k.b) $$
o\`u \ $S_j \in \C[[b]]$ \ v\'erifie \ $S_j(0) = 1$, avec la convention \ $\alpha_j = \emptyset$ \ si \ $p_j = 0$.\\
Rappelons (voir [B.09b])  qu'un th\`eme \ $[\lambda]-$primitif  de rang k est toujours isomorphe \`a un th\`eme du type ci-dessus, les \ $S_j$ \ pouvant \^etre choisis polynomiaux en \ $b$ \ (avec \ $S_j(0) = 1$ ) \  et on peut borner les degr\'es des \ $S_j$ \ en fonction des invariants fondamentaux (voir les familles standards dans [B.09b]).

\section{Le th\'eor\`eme du changement de variable.}

\begin{lemma}\label{chgt 1}
Soit \ $\theta \in \C[[a]]$ \ v\'erifiant  \ $\theta(0) = 0$ \ et \ $\theta'(0) \not= 0 $. Alors on a un  (unique) homomorphisme de \ $\C-$alg\`ebre unitaire 
$$ \Theta : \hat{A} \to \hat{A} $$
v\'erifiant \ $\Theta(a) : = \alpha = \theta(a)$ \ et \ $\Theta(b) : = \beta = b.\theta'(a) $.
\end{lemma}

\parag{Preuve} Les relations \ $a^n.b = b.a^n + n.b.a^{n-1}.b$ \ valables pour tout \ $n \in \mathbb{N}$ \ donnent la relation \ $\theta(a).b = b.\theta(a) + b.\theta'(a).b$ \ dans \ $\hat{A}$ \ pour toute s\'erie formelle \ $\theta \in \C[[a]]$. Alors on aura, toujours dans \ $\hat{A}$ :
  $$\alpha.\beta = \theta(a).b.\theta'(a) = (b.\theta(a) + b.\theta'(a).b).\theta'(a) = \beta.\alpha + \beta^2 $$
  ce qui montre notre assertion. $\hfill \blacksquare$\\

Nous appellerons {\bf changement de variable} un \ $\theta \in \C[[a]]$ \ v\'erifiant  \ $\theta(0) = 0$ \ et \ $\theta'(0) \not= 0 $. 

\begin{defn}\label{chgt 2}
Soit \ $E$ \ un (a,b)-module r\'egulier (ou m\^eme,  plus g\'en\'eralement,  un \ $\hat{A}-$module). On d\'efinit, pour tout changement de variable \ $\theta$, un nouvel (a,b)-module \ $\theta_*(E)$, qui sera appel\'e le {\bf chang\'e de variable de \ $E$ \ par \ $\theta$},  comme \'etant le (a,b)-module (resp. le \ $\hat{A}-$module) obtenu en faisant agir sur \ $E$ \ l'\'el\'ement \ $a \in \tilde{A}$ \ par \ $\alpha : = \theta(a)$ \ et l'\'el\'ement \ $b \in \tilde{A}$ \ par \ $\beta : = b.\theta'(a)$.
\end{defn}

Ceci revient \`a faire agir \ $\hat{A}$ \ sur \ $E$ \ via l'automorphisme \ $\Theta$, c'est-\`a-dire que pour \ $u \in \hat{A}$ \ et \ $x \in E$ \ on pose \ $u._{\theta}x : = \Theta(u).x$.\\

\begin{prop}\label{chgt 3}
Soit \ $E$ \ un (a,b)-module \`a p\^ole simple, alors \ $\theta_*(E)$ \ est \`a p\^ole simple, on a \ $b^n.E = \beta^n.E = b^n.\theta_*(E)$ \ pour chaque entier \ $n$,  et l'action de \ $b^{-1}.a$ \ sur \ $\theta_*(E)\big/b.\theta_*(E)$ \ est la m\^eme\footnote{c'est en fait l'action de \ $\beta^{-1}.\alpha$ \ sur l'espace vectoriel \ $E\big/\beta.E \simeq E\big/b.E$.} que celle de \ $b^{-1}.a$ \ sur \ $E\big/b.E$. \\
Soit \ $E$ \ un (a,b) r\'egulier  et soit \ $\theta$ \ un changement de variable. Alors \ $\theta_*(E)$ \ est un (a,b)-module r\'egulier de m\^eme rang. Ses sous-(a,b)-modules sont les m\^emes que ceux de \ $E$, et ils sont normaux dans \ $\theta_*(E)$ \ si et seulement s'ils le sont dans \ $E$. Le satur\'e par \ $b^{-1}.a$ \ de \ $\theta_*(E)$ \ est \ $\theta_*(E^{\sharp})$, o\`u \ $E^{\sharp}$ \ est le satur\'e de \ $E$ \ par \ $b^{-1}.a$ \ et \ $\theta_*(E)$ \  a le m\^eme polyn\^ome de Bernstein que \ $E$.
\end{prop}

\parag{Preuve} Si \ $F \subset E$ \ est stable par \ $a$ \ et \ $b$, il est r\'egulier car \ $E$ \ l'est, et il est donc stable par \ $\alpha$ \ et \ $\beta$. La r\'eciproque r\'esulte du fait qu'il existe un changement de variable \ $\eta \in \C[[\alpha]]$ \ tel que \ $a = \eta(\alpha)$.\\
Si \ $F$ \ est normal dans \ $E$, soit \ $x \in E$ \ v\'erifiant \ $\beta.x \in F$. On a donc \ $\theta'(a).x \in F$, et comme \ $\theta'(a)$ \ est un inversible de \ $\C[[a]]$, on en conclut que \ $\theta_*(F)$ \ est normal dans \ $\theta_*(E).$\\
Montrons que \ $\theta_*(E)$ \ est \`a p\^ole simple si \ $E$ \ l'est. Mais comme \ $\theta'(a)$ \ est inversible dans \ $\hat{A}$ \ et \ $E$ \ r\'egulier, on aura \ $\beta.E = b.E$. Alors \ $\alpha.E \subset a.E \subset b.E = \beta.E$, et \ $\theta_*(E)$ \ est bien \`a p\^ole simple.\\
Comme l'espace vectoriel \ $E\big/b.E$ \ est invariant par changement de variable, montrons  l'invariance de l'action de \ $b^{-1}.a$ \ sur cet espace vectoriel. On se ram\`ene imm\'ediatement au cas \ $\theta'(0) = 1$. On a alors en posant \ $\theta(a) = a + a^2.\eta(a)$
 $$\beta^{-1}.\alpha = \theta'(a)^{-1}.b^{-1}.\theta(a) = (1 + a^2.\eta'(a) + 2a.\eta(a) )^{-1}.(b^{-1}.a).(1 + a.\eta(a)) \quad {\rm dans} \quad E $$
 et comme \ $a$ \ induit l'application nulle sur \ $E\big/b.E$, l'invariance en d\'ecoule.\\
Si \ $E$ \ est r\'egulier, il existe une inclusion de \ $E$ \ dans un (a,b)-module \ $G$ \ \`a p\^ole simple. On a alors \ $\theta_*(E) \subset \theta_*(G)$ \ qui est \`a p\^ole simple. D'o\`u la r\'egularit\'e de \ $\theta_*(E)$. On a obtenu de plus que le satur\'e de \ $\theta_*(E)$ \ est contenu dans \ $\theta_*(E^{\sharp})$, o\`u \ $E^{\sharp}$ \ d\'esigne le satur\'e de \ $E$ \ par \ $b^{-1}.a$ \ dans  $E[b^{-1}]$. Pour obtenir l'\'egalit\'e consid\'erons un sous-(a,b)-module \`a p\^ole simple \ $F$ \ v\'erifiant  \ $\theta_*(E) \subset F \subset \theta_*(E^{\sharp})$. Alors le changement de variable inverse \ $\theta^{-1}$ \ donnera \ $E \subset \theta^{-1}_*(F) \subset E^{\sharp}$. Comme \ $\theta^{-1}_*(F)$ \ est \`a p\^ole simple, on aura, puisque \ $E^{\sharp}$ \ est le plus petit module \`a p\^ole simple contenant \ $E$,  \ $\theta^{-1}_*(F) = E^{\sharp}$, d'o\`u l'\'egalit\'e \ $\theta_*(E^{\sharp}) = \theta_*(E)^{\sharp}$.$\hfill \blacksquare$\\

\parag{Remarques}
\begin{enumerate}
\item On prendra garde que si \ $F \subset E$ \ est un sous-(a,b)-module, il est stable par \ $\alpha$ \ et \ $\beta$, mais le sous-(a,b)-module ainsi obtenu (qui est \ $\theta_*(F)$) \ n'est pas, \`a priori, isomorphe \`a \ $F$. Le lemme suivant montrera que pour les (a,b)-modules r\'eguliers de rang 1, \ $\theta_*(E)$ \ est toujours isomorphe \`a \ $E$. On verra plus loin que c'est aussi le cas pour th\`eme \ $[\lambda]-$primitif de rang 2, mais que c'est en g\'en\'eral faux pour un th\`eme \ $[\lambda]-$primitif de rang \ $\geq 3$.
\item  Comme \ $a.E + b.E$ \ est un sous-(a,b)-module de \ $E$ \ il est imm\'ediat de voir que \ $a.\theta_*(E) + b.\theta_*(E) = \theta_*(a.E + b.E)$. Un (a,b)-module r\'egulier \ $E$ \ est monog\`ene si et seulement si \ $\dim_{\C}(E\big/(a.E + b.E)) = 1$. Cette propri\'et\'e \'etant invariante par changement de variable, on en conclut que si \ $E$ \ est r\'egulier monog\`ene, pour tout changement de variable \ $\theta$ \ le (a,b)-module \ $\theta_*(E)$ \ est encore r\'egulier monog\`ene. Donc \ $\theta_*(E)$ \ a le m\^eme \'el\'ement de Bernstein que \ $E$.\\
 On notera que ceci est une cons\'equence simple du premier th\'eor\`eme de structure des (a,b)-modules monog\`enes r\'eguliers de [B. 09a]  puisque pour \ $u \in \hat{A}$ \ la forme initiale en (a,b) de \ $\Theta(u)$ \ est la m\^eme que celle de \ $u$, \`a un scalaire non nul pr\`es (qui est \ $\theta'(0)^k$ \ si la forme initiale est de degr\'e \ $k$).
 \end{enumerate}

\begin{lemma}\label{chgt 4}
Soit \ $\lambda \in \C$ \ et \ $E_{\lambda} : = \C[[b]].e_{\lambda}$ \ le (a,b)-module de rang 1 d\'efini par la relation \ $ a.e_{\lambda} = \lambda.b.e_{\lambda}$\footnote{donc \ $E_{\lambda} \simeq \tilde{A}\big/ \tilde{A}.(a - \lambda.b) \simeq \hat{A}\big/\hat{A}.(a - \lambda.b).$}. Pour tout changement de variable \ $\theta$ \ on a \ $\theta_*(E_{\lambda}) \simeq E_{\lambda}$.
\end{lemma}

\parag{Remarque} Comme on a \ $Aut_{a,b}(E_{\lambda}) = \C^*$, c'est-\`a-dire que le groupe des automorphismes du (a,b)-module \ $E_{\lambda}$ \ est r\'eduit aux homoth\'eties non nulles, l'isomorphisme entre \ $\theta_*(E_{\lambda})$ \ et \ $E_{\lambda}$ \ est unique \`a un scalaire multiplicatif non nul pr\`es.

\parag{Preuve} On a d\'ej\`a vu que si \ $E$ \ est \`a p\^ole simple \ $\theta_*(E)$ \ est \'egalement \`a p\^ole simple et il est de m\^eme rang que \ $E$ \ car \ $\beta.E = b.E$ \ montre que l'espace vectoriel \ $E\big/b.E$ \ est invariant par changement de variable. Donc \ $\theta_*(E_{\lambda})$ \ est isomorphe \`a \ $E_{\mu}$ \ pour un \ $\mu \in \C$ \ d'apr\`es la classification des (a,b)-modules r\'eguliers de rang 1. Mais l'\'egalit\'e des polyn\^omes de Bernstein donne\ $\mu = \lambda$. $\hfill \blacksquare$

\parag{Remarque} Une cons\'equence simple de ce qui pr\'ec\`ede est que si \ $(F_j)_{j \in [0,k]}$ \ est une suite de Jordan-H{\"o}lder d'un (a,b)-module r\'egulier \ $E$ \ de rang \ $k$, les \ $\theta_*(F_j)_{j \in [0,k]}$ \ forment une suite de Jordan-H{\"o}lder de \ $\theta_*(E)$ \ et les quotients \ $\theta_*(F_{j+1})\big/\theta_*(F_j)$ \ sont isomorphes aux quotients  \ $F_{j+1}\big/F_j$.\\

\begin{lemma}\label{Xi}
 Soit \ $\lambda \in \C\setminus -\mathbb{N}$.  Il existe un unique (a,b)-module \ $E$ \ \`a p\^ole simple de rang 2  tel que l'action de \ $a$ \ dans \ $E\big/b^2.E$ \ soit \'egale \`a \ $b.\left(\begin{matrix} \lambda & 1\\ 0 & \lambda\end{matrix}\right)$. Il est invariant par tout changement de variable et c'est le (a,b)-module  $$\Xi_{\lambda}^{(1)}\simeq \C[[b]].s^{\lambda-1}.Log\, s \oplus \C[[b]].s^{\lambda-1}.$$
 \end{lemma}

\parag{Preuve} On sait que \ $E$ \ admet un sous-(a,b)-module normal isomorphe \`a \ $E_{\lambda}$, et que, quitte \`a modifier le second vecteur de base modulo \ $b^2.E$ \ on peut supposer que \ $a.e_2 = \lambda.be_2$\footnote{gr\^ace \`a la proposition 1.3  de [B.93].}; posons  alors \ $a.e_1 = \lambda.b.e_1 + b.e_2 + b^2.S(b).e_1 + b^2.T(b).e_2 $. Si on cherche une base \ $\varepsilon_1, \varepsilon_2$ \ sur \ $\C[[b]]$ \ sous la forme
\begin{equation*}
 \varepsilon_1 = e_1 + b.U.e_1 + b.V.e_2 \quad {\rm et}\quad  \varepsilon_2 = e_2  \tag{*}
 \end{equation*}
v\'erifiant \ $a.\varepsilon_1 = \lambda.b.\varepsilon_1 + b.\varepsilon_2$ \ et \ $a.\varepsilon_2 = \lambda.b\varepsilon_2$ \ (seconde condition que l'on a d\'ej\`a suppos\'ee v\'erifi\'ee), on trouve pour \ $U \in \C[[b]]$ \ l'\'equation diff\'erentielle
\begin{equation*}
 (1 + b.S).U + b.U' = -S \tag{A}
 \end{equation*}
puis, $U$ \ \'etant ainsi choisi, l'\'equation diff\'erentielle
\begin{equation*}
 b.V' + V = - U.(1 + b.T) - T . \tag{B}
 \end{equation*}

Explicitons les calculs :  Pour r\'esoudre \ $(A)$, notons \ $\Sigma$ \ la primitive sans terme constant de \ $S$, et posons \ $U : = \Gamma.exp(-\Sigma)$. Alors \ $(A)$ \ devient 
\begin{equation*}
\Gamma + b.\Gamma' = -S.exp(\Sigma) \tag{A'}
\end{equation*}
qui a une unique solution dans \ $\C[[b]]$.\\
La r\'esolution de \ $(B)$ \ est alors facile.
On en conclut que tout tel (a,b)-module est isomorphe \`a celui d\'efini par la \ $\C[[b]]-$base \ $\varepsilon_1, \varepsilon_2$ \ et les relations \ $(^*)$. Pour \ $\lambda \not\in - \mathbb{N}$ \  c'est bien s\^ur \ $\Xi_{\lambda}^{(1)}$ \ avec \ $\varepsilon_1 : = s^{\lambda-1}.Log\, s$ \ et \ $\varepsilon_2 : = s^{\lambda-1}$. \\
Comme un changement de variable  de la forme \ $\theta(a) = a + a^2.\eta(a)$ \ ne change pas l'action de \ $a$ \ et \ $b$ \ sur \ $E\big/b^2.E$ \ quand \ $E$ \ est \`a p\^ole simple, l'invariance par ces changements de variables r\'esulte de la caract\'erisation pr\'ec\'edente. Pour les changement de variables de la forme \ $\theta(a) = \rho.a$ \ avec \ $\rho \in \C^*$, l'assertion est imm\'ediate.$\hfill \blacksquare$

\begin{lemma}\label{chgt 7}
Soit \ $E$ \ un (a,b)-module monog\`ene r\'egulier \ $[\lambda]-$primitif de rang \ $k \geq 3$ \ et soit \ $(F_j)_{j\in [1,k]}$ \ une suite de Jordan-H{\"o}lder de \ $E$. Supposons que \ $F_{k-1}$ \ soit un th\`eme ainsi que \ $E\big/F_{k-2}$. Alors \ $E$ \ est un th\`eme.
\end{lemma}

\parag{Preuve}
Soit \ $E_{\lambda} \subset E$ \ un sous-(a,b)-module normal de rang 1. D'apr\`es le th\'eor\`eme 2.1.6 de [B.09b],  il suffit de montrer que \ $E_{\lambda} = F_1$ \ pour conclure. Consid\'erons l'image de \ $E_{\lambda}$ \ dans le quotient \ $E\big/F_{k-2}$. Le normalis\'e de cette image est de rang \ $\leq 1$. Ceci montre que cette image est contenue dans \ $F_{k-1}\big/F_{k-2}$ \ qui est l'unique sous-module normal de rang 1 de \ $E\big/F_{k-2}$ \ puisque l'on a suppos\'e que c'est un th\`eme. Alors on a \ $E_{\lambda} \subset F_{k-1}$ \ qui est un th\`eme. Donc \ $E_{\lambda} = F_1$ \ qui est l'unique sous-module normal de rang 1  de ce th\`eme. $\hfill \blacksquare$\\

\begin{cor}\label{chgt 5}
Soit \ $E$ \ un th\`eme \ $[\lambda]-$primitif de rang k et soit \ $\theta$ \ un changement de variable. Alors \ $\theta_*(E)$ \ est un th\`eme \ $[\lambda]-$primitif de rang k. Il a m\^emes invariants fondamentaux que \ $E$.
\end{cor}

\parag{Preuve} En rang 1 c'est clair. En rang 2 \'egalement car si \ $j : E \to \Xi^{(1)}_{\lambda}$ \ est une injection (a,b)-lin\'eaire, alors \ $j : \theta_*(E) \to \theta_*(\Xi_{\lambda}^{(1)})$ \ sera une injection \ (a,b)-lin\'eaire, et on a \ $\theta_*(\Xi_{\lambda}^{(1)}) \simeq \Xi_{\lambda}^{(1)}$ \ d'apr\`es le lemme pr\'ec\'edent.\\
Montrons maintenant le r\'esultat par r\'ecurrence sur \ $k \geq 2$. Supposons le r\'esultat d\'emontr\'e pour le rang \ $k \geq 2$ \ et montrons-le en rang \ $k +1$. Soit \ $E$ \ un th\`eme \ $[\lambda]-$primitif de rang \ $k+1$. Notons \ $F_j$ \ le sous-(a,b)-module normal de rang \ $j$ \ de \ $E$. Alors \ $F_k$ \ est un th\`eme \ $[\lambda]-$primitif de rang \ $k$ \ et donc \'egalement \ $\theta_*(F_k)$. Comme \ $E\big/F_{k-1}$ \ est un th\`eme \ $[\lambda]-$primitif de rang \ $2$, il en est de m\^eme de \ $\theta_*(E)\big/\theta_*(F_{k-1}) \simeq \theta_*(E\big/F_{k-1})$. Comme \ $\theta_*(E)$ \ est monog\`ene \ $[\lambda]-$primitif et admet la suite de Jordan-H{\"o}lder \ $\theta_*(F_j), j \in [1,k+1]$, le lemme suivant permet de conclure que \ $\theta_*(E)$ \ est un th\`eme.\\
Les invariants fondamentaux \'etant donn\'es par les quotients succesifs de l'unique suite de Joran-H{\"o}lder de \ $\theta_*(E)$, la conclusion r\'eulte de la remarque qui pr\'ec\`ede le lemme \ref{Xi}. \  $\hfill \blacksquare$\\

\begin{defn} Soit \ $f : E \to F$ \ une application   \ $\C-$lin\'eaire entre deux  \ $\C[[b]]-$modules. On dira que \ $f$ \ est {\bf \ $b-$compatible} si on a pour chaque entier \ $n$ \ l'inclusion \ $f(b^n.E) \subset b^n.F$. Si \ $f$ \ est bijective, on dira que \ $f$ \ est {\bf strictement $b-$compatible} si on a \ $f(b^n.E) = b^n.F$ \ pour tout \ $n$. Ceci \'equivaut \`a demander que \ $f$ \ et \ $f^{-1}$ \ soient \ $b-$compatibles.
\end{defn}

\parag{Remarque} Il n'existe pas\footnote{On voit facilement que la matrice est triangulaire dans la base \ $b^n, n \geq 0$.} d'application \ $b-$compatible bijective de \ $\C[[b]]$ \ dans \ $\C[[b]]$ \ qui ne soit pas strictement compatible.

\begin{defn} Soit \ $E$ \ un \ $\C[[b]]-$module libre de type fini et soit \ $f_{\theta} : E \to E$ \ une application \ $\C-$lin\'eaire d\'ependant d'un param\`etre \ $\theta \in \C^N$. On suppose que\ $f$ \ est \ $b-$compatible pour chaque valeur de \ $\theta$. On dira que \ $f$ \ {\bf d\'epend  polynomialement de  \ $\theta$} \ si pour chaque \ $x \in E$ \ et chaque  entier \ $n$ \ l'application \ $\C^N \to E\big/b^n.E$ \ qui \`a \ $\theta$ \ associe la classe de \ $f_{\theta}(x)$ \ dans \ $E\big/b^n.E$ \ est polynomiale.
\end{defn}

\parag{Remarques}
\begin{enumerate}
\item Dans le cas d'un \'el\'ement \ $S_{\theta}$ \ de \ $\C[[b]]$ \ d\'ependant du param\`etre \ $\theta \in \C^N$ \ il est \'equivalent de demander que chaque coefficient de \ $S_{\theta}$ \ soit un polyn\^ome en \  $\theta$ \ ou que l'op\'erateur de multiplication par \ $S_{\theta}$ \ de \ $\C[[b]]$ \ dans \ $\C[[b]]$ \ soit \ $b-$compatible.
\item On notera que si l'on a deux applications \ $\C-$lin\'eaires b-compatibles  \ $f$ \ et \ $g$ \ de \ $\C[[b]]$ \ dans \ $\C[[b]]$, d\'ependant polynomialement  d'un param\`etre \ $\theta \in \C^N$, leur compos\'ee est encore lin\'eaire \ $b-$compatible et d\'epend polynomialement de \ $\theta$. En effet le coefficient de \ $b^n$ \ de \ $g( f(x))$ \ ne d\'epend que des coefficients de \ $f(x)$ \ dans \ $\C[[b]]\big/b^{n+1}.\C[[b]]$ \ et de l'endomorphisme induit par \ $g$ \ sur ce quotient dont la matrice est \`a coefficients dans les polyn\^omes en \ $\theta$.
 \end{enumerate}

Dans la suite on consid\`erera des changements de variable de la forme
$$ \theta(a) = a + \sum_{j=2}^N \ \theta_j.a^j $$
o\`u l'on prendra  \ $(\theta_2, \dots, \theta_N) \in \C^{N-1}$ \  comme  param\`etre. On va d\'ej\`a \'etudier ces changements de variables (donc \ $\alpha_{\theta} = \theta(a)$ \ et \ $\beta_{\theta} = b.\theta'(a)$) \  sur le (a,b)-module \ $E_{\lambda}$ \ pour \ $\lambda \in \C \setminus -\mathbb{N}$.

\begin{lemma}
 Soit \ $E : = E_{\lambda} = \C[[b]].e_{\lambda}$ \ o\`u \ $a.e_{\lambda} = \lambda.b.e_{\lambda}$. Il existe pour chaque \ $n \in \mathbb{N}$ \ un \'el\'ement \ $\chi_n \in \C[[u]]$ \ d\'ependant polynomialement de \ $\theta \in \C^{N-1}$ \ tel que l'on ait dans \ $E_{\lambda}$
\begin{equation*}
 b^n.e_{\lambda} = \beta^n.\chi_n(\beta).e_{\lambda} \tag{1}
 \end{equation*}
o\`u l'operateur \ $\beta : = \beta_{\theta}$ \ est d\'efini par \ $\beta : = b.\theta'(a)$ \ sur \ $E_{\lambda}$.
\end{lemma}

\parag{Preuve} Notons d\'ej\`a que \ $\alpha = \theta(a)$ \ est une application \ $\C-$lin\'eaire $b-$compatible d\'ependant polynomialement\footnote{et m\^eme de fa{\c c}on affine de \ $(\theta_2, \dots, \theta_N)$.} de \ $\theta\in \C^{N-1}$ \ sur tout (a,b)-module \`a p\^ole simple. Il en est alors de m\^eme pour \ $\beta$ \ qui est m\^eme la compos\'ee de \ $b$ \ et de l'application lin\'eaire (bijective) strictement compatible \ $\theta'(a)$ \ qui d\'epend polynomialement de \ $\theta$. Comme l'\'egalit\'e \ $\C[[b]].e_{\lambda} = \C[[\beta]].e_{\lambda}$ \ et la $b-$compatibilit\'e assurent de l'existence et de l'unicit\'e de \ $\chi_n$ \ pour tout entier \ $n$, il nous suffit de montrer que \ $\chi_n$ \ d\'epend polynomialement de \ $\theta$. Fixons \ $p \in \mathbb{N}, p \gg 1$. Alors on a deux bases de l'espace vectoriel  \ $V_p : = E_{\lambda}\big/b^p.E_{\lambda}$, la base \ $b^q.e_{\lambda}$ \ et la base \ $\beta^q.e_{\lambda}, q \in [0,p-1]$. Et on a clairement \ $\beta^q.e_{\lambda} = b^q.e_{\lambda} + b^{q+1}.V$ \ pour tout \ $q$. Ceci montre que ce changement de base est triangulaire avec des \ $1$ \ sur la diagonale. Donc de d\'eterminant \'egal \`a \ $1$.\\
Par ailleurs, comme \ $\beta$ \ d\'epend polynomialement de \ $\theta$, la matrice de ce changement de base est polynomiale en \ $\theta$. Sa matrice inverse \'egalement puisque le d\'eterminant vaut \ $1$. Ceci permet de conclure, car pour calculer le coefficient de \ $\beta^m$ \ dans  \ $\chi_n$, on peut se contenter de travailler dans \ $V_{n+m+1}$.$\hfill \blacksquare$

\begin{cor}
Dans la situation pr\'ec\'edente il existe une unique \ $S^{\lambda}_{\theta} \in \C[[b]]$ \ v\'erifiant \ $S^{\lambda}_{\theta}(0) =1$ \ et d\'ependant polynomialement de \ $\theta \in \C^{N-1}$ \ tel que l'on ait $$\alpha_{\theta}.\varepsilon_{\lambda}^{\theta} = \lambda.\beta_{\theta}.\varepsilon_{\lambda}^{\theta} $$
si l'on pose \ $\varepsilon_{\lambda}^{\theta} = S^{\lambda}_{\theta}(\beta_{\theta}).e_{\lambda}$.
\end{cor}

\parag{Preuve} Comme on suppose \ $\lambda$ \ fix\'e dans la suite, on omettera l'exposant \ $\lambda$ \ pour les \'el\'ements de \ $\C[[\beta]]$ \ que l'on consid\`erera. Commen{\c c}ons par remarquer que l'on a
$$ \alpha.e_{\lambda} = \lambda.b.e_{\lambda} + \sum_{j=2}^{N} \theta_j.\lambda_j.b^j.e_{\lambda} $$
o\`u \ $\lambda_j = \lambda.(\lambda+1)\dots (\lambda+j-1) $ \ pour \ $j \geq 1$, ce qui, en utilisant le lemme pr\'ec\'edent nous fournit \ $R_{\theta} \in \C[[\beta]]$ \  d\'ependant polynomialement de \ $\theta \in \C^{N-1}$ \ tel que
$$ \alpha.e_{\lambda} = \lambda.\beta.e_{\lambda} + \beta^2.R_{\theta}(\beta).e_{\lambda}. $$
L'\'egalit\'e
$$ \alpha.S_{\theta}.e_{\lambda} = \lambda.\beta.S_{\theta}.e_{\lambda} $$
donne alors
\begin{equation*}
 S_{\theta}.\big[ \lambda.\beta.e_{\lambda} + \beta^2.R_{\theta}(\beta).e_{\lambda}\big] + \beta^2.S_{\theta}'.e_{\lambda} =  \lambda.\beta.S_{\theta}.e_{\lambda} 
 \end{equation*}
 ce qui, apr\`es simplification donne l'\'equation diff\'erentielle
 $$ S_{\theta}' + R_{\theta}.S_{\theta} = 0 $$
et donc  \ $S_{\theta} = exp\big[ - \tilde{R}_{\theta} \big] $ \ o\`u \ $\tilde{R}_{\theta}$ \ d\'esigne la primitive nulle en \ $0$ \ de \ $R_{\theta}$. Comme \ $R_{\theta}$ \ d\'epend polynomialement de \ $\theta$ \ il en sera de m\^eme de \ $\tilde{R}_{\theta}$ \ et aussi de \ $exp\big[ - \tilde{R}_{\theta} \big] $ \ puisque la nullit\'e du terme constant de \ $\tilde{R}_{\theta}$ \ assure que le coefficient de \ $\beta^p$ \ dans l'exponentielle ne d\'epend que du d\'eveloppement limit\'e  \`a l'ordre \ $p$ \ de l'exponentielle. $\hfill \blacksquare$.
 
 \begin{prop}\label{crucial}
 Soit \ $E : = \tilde{A}\big/\tilde{A}.(a - \lambda.b)(1 + z.b^p)^{-1}.(a - (\lambda+p-1).b) $ \ o\`u l'on suppose \ $\lambda \not\in - \mathbb{N}$ \ et \ $p \in \mathbb{N}^*$. On effectue sur \ $E$ le changement de variable \ $\theta(a) = a + \sum_{j=2}^N \theta_j.a^j $. Alors \ $\theta_*(E)$ \ est isomorphe au (a,b)-module 
  $$\tilde{A}\big/\tilde{A}.(a - \lambda.b)(1 + z_{\theta}.b^p)^{-1}.(a - (\lambda+p-1).b),$$
   o\`u \ $z_{\theta}$ \ est un polyn\^ome en \ $\theta \in \C^{N-1}$.
 \end{prop}
 
 \parag{Remarque} Si on part de \ $z \not= 0$ \ et \ $\lambda_1 \in 1 + \mathbb{Q}^{+*}$, alors \ $E$ \ est un th\`eme de rang 2. Si l'on sait par ailleurs que pour tout \ $\theta$ \ le (a,b)-module \ $\theta_*(E)$ \ est encore un th\`eme, cela impliquera que \ $z_{\theta}$ \ n'est pas nul. Or les seuls polyn\^omes qui ne s'annulent jamais sur \ $\C^{N-1}$ \ sont les polyn\^omes constants. Cela montrera que l'on a \ $z_{\theta} \equiv z$ \ pour tout \ $\theta$.
 Ce raisonnement sera la clef de la d\'emonstration du th\'eor\`eme qui suit.\\

\begin{thm}\label{th. fond.}
Soit \ $E$ \ un th\`eme \ $[\lambda]-$primitif de rang \ $k$ \ et d'invariants fondamentaux \ $\lambda_1, p_1, \dots, p_{k-1}$. Notons \ $z_1, \dots, z_{k-1}$ \ les param\`etres des th\`emes  \ $[\lambda]-$primitif de rang 2 \ $F_2, F_3\big/F_1, \dots, F_k\big/F_{k-2}$. Soit \ $\theta$ \ un changement de variable de la forme \ $\theta(a) = r.a + \sum_{j=2}^{\infty} \theta_j.a^j$ \ avec \ $r \not= 0$. Alors \ $\theta_*(E)$ \ est un th\`eme \ $[\lambda]-$primitif de rang \ $k$ \ et d'invariants fondamentaux \ $\lambda_1, p_1, \dots, p_{k-1}$. Les param\`etres des th\`emes  \ $[\lambda]-$primitif de rang 2 
 $$\theta_*(F_2), \theta_*(F_3\big/F_1), \dots, \theta_*(F_k\big/F_{k-2})$$
 sont les nombres \ $r^{p_1}.z_1, \dots, r^{p_{k-1}}.z_{k-1}$.
\end{thm}

\parag{D\'emonstration} Commen{\c c}ons par montrer qu'il suffit de prouver le r\'esultat en rang 2. En effet si l'on sait d\'ej\`a que si \ $E$ \ est un th\`eme  \ $[\lambda]-$primitif  de rang 2 d'invariants fondamentaux \ $\lambda_1, p_1$ \  et de param\`etre \ $z$, alors \ $\theta_*(E)$ \ est un th\`eme  \ $[\lambda]-$primitif de rang \ $2$ \  ayant m\^emes invariants fondamentaux et de param\`etre \ $r^{p_1}.z$, on en d\'eduit les propri\'et\'es suivantes  pour un th\`eme \ $[\lambda]-$primitif de rang \ $k$ \ et d'invariants fondamentaux \ $\lambda_1, p_1, \dots, p_{k-1}$ :
\begin{enumerate}[i)]
\item D'abord \ $\theta_*(E)$ \ est monog\`ene r\'egulier de rang \ $k$. Si on suppose que pour un th\`eme \ $[\lambda]-$primitif  \ $F$ \ de rang \ $k-1$ \ $\theta_*(F)$ \ est encore un th\`eme \ $[\lambda]-$primitif  \ $F$ \ de rang \ $k-1$ \ avec les m\^emes invariants fondamentaux que \ $F$, alors le lemme \ref{chgt 7} permet de montrer qu'il en sera de m\^eme en rang \ $k$. Le rang 1 \'etant clair, il suffit donc bien de traiter le cas du rang 2.
\item De m\^eme la quasi-invariance\footnote{par rapport au caract\`ere \ $\theta \mapsto \theta'(0)^p $ \ du groupe des changements de variables,  pour un th\`eme de rang \ $2$ \ d'invariants fondamentaux \ $(\lambda_1, p)$.}des param\`etres des th\`emes  \ quotients \ $F_{j+2}\big/F_j$ \ de rang \ $2$ \ est cons\'equence de ce r\'esultat  pour \ $k = 2$, puisque \ $\theta_*(F_j)$ \ est le sous-th\`eme normal de rang \ $j$ \ de \ $\theta_*(E)$.
\end{enumerate}

Pour traiter le cas de rang 2, montrons d\'ej\`a que si \ $E$ \ est un th\`eme \ $[\lambda]-$primitif  de rang $2$ \  et si \ $\theta$ \ est un changement de variable, alors \ $\theta_*(E)$ \ est un th\`eme  \ $[\lambda]-$primitif de rang \ $2$. On a une injection (a,b)-lin\'eaire \ $j : E \to \Xi_{\lambda}^{(1)}$. On aura donc une injection \ $j : \theta_*(E) \to \theta_*(\Xi_{\lambda}^{(1)})$. Comme le lemme \ref{Xi} donne un isomorphisme \ $\theta_*(\Xi_{\lambda}^{(1)}) \simeq \Xi_{\lambda}^{(1)}$, cela montre que \ $\theta_*(E)$ \ est un th\`eme \ $[\lambda]-$primitif.\\
La quasi-invariance du param\`etre est cons\'equence alors du fait que, pour un changement de variable polynomial\footnote{c'est-\`a-dire de la forme \ $\theta(a) = r.a + \sum_{j=2}^N \theta_j.a^j$.}, le param\`etre est une fonction polynomiale de \ $\theta$ \ qui ne peut s'annuler (sinon on n'aurait pas un th\`eme !). Ceci implique aussi le cas d'un changement de variable g\'en\'eral, car pour \ $\theta$ \ donn\'e et \ $E$ \ donn\'e, en rempla{\c c}ant \ $\theta$ \ par \ $\theta_N$ \ son d\'eveloppement limit\'e d'ordre \ $N \gg1$, on aura que \ $(\theta_N)_*(E)$ \ est isomorphe \`a \ $\theta_*(E)$. D'o\`u le cas g\'en\'eral. $\hfill \blacksquare$

 \parag{Exemple} Soit \ $\varphi : = s^{\lambda+p-2}.Log\, s + S(b).s^{\lambda-2} $ \ o\`u \ $\lambda \in 1+ \mathbb{Q}^{*+},  S \in \C[[b]], S(0) \not= 0$ \ et \ $p \in \mathbb{N}^*$. Alors \ $E : = \A.\varphi \subset \Xi_{\lambda}^{(1)}$ \ est un th\`eme de rang 2. Nous allons calculer pour \ $p \geq 1$ \ le param\`etre de ce th\`eme dont les invariants fondamentaux sont \ $\lambda_1 = \lambda, \lambda_2 = \lambda+p-1$ \ et \ $p_1 = p$.\\
 On a 
 \begin{align*}
& (a - (\lambda+p-1).b).\varphi = \frac{s^{\lambda+p-1}}{\lambda+p-1} + S(b).s^{\lambda -1} + b^2.S'(b).s^{\lambda-2} - (\lambda+p-1).S(b).\frac{s^{\lambda-1}}{\lambda -1} \\
&\qquad\qquad = \big[(\lambda-1).S(b) + b.S'(b) - (\lambda+p-1).S(b)\big].\frac{s^{\lambda -1}}{\lambda-1}  + \rho.b^p.\frac{s^{\lambda-1}}{\lambda - 1} 
\end{align*}
o\`u l'on a pos\'e \ $\rho = \frac{\Gamma(\lambda+p-1)}{\Gamma(\lambda-1)}$ \ puisque \ $$b^p.\frac{s^{\lambda-1}}{\lambda-1}  = \frac{\Gamma(\lambda-1)}{\Gamma(\lambda+p)}.s^{\lambda+p-1}.$$
On arrive donc, pour \ $p \geq 1$ \ \`a l'\'egalit\'e
$$ (a - (\lambda+p-1).b).\varphi = \big[b.S'(b) - pS(b) + \rho.b^p \big].\frac{s^{\lambda-1}}{\lambda-1} $$
ce qui montre que \ $ (a - \lambda.b).T^{-1}(b).(a - (\lambda+p-1).b).\varphi = 0$ \ o\`u l'on a pos\'e \ $T(b) : = (b.S'(b) - pS(b) + \rho.b^p)$. On notera que comme on a suppos\'e \ $S(0) \not= 0$ \ et \ $p \geq 1$ \ on a \ $T(0) = -p.S(0) \not= 0 $.\\
Le param\`etre fondamental de ce th\`eme de rang 2 est donc \'egal \`a  
$$\alpha : = -\frac{\rho}{p.S(0)} = -\frac{\Gamma(\lambda+p-1)}{p.\Gamma(\lambda-1).S(0)} $$ 
 puisque la s\'erie formelle \ $b.S'(b) - p.S(p)$ \ a un coefficient de \ $b^p$ \ qui est nul.

\section{Un contre-exemple.}

 Nous allons montr\'e sur un exemple de rang 3  que les param\`etres non principaux ne sont pas. en g\'en\'eral, invariants par changements de variables.
 
 Commen{\c c}ons par montrer une proposition.
 
 \begin{prop} Soit \ $\lambda \in \C$ \ tel que \ $\lambda \in  2+\mathbb{Q}^{*+} $. Consid\'erons alors l'\'el\'ement de \ $\Xi_{[\lambda]}^{(2)}$
  $$ e : =  s^{\lambda-1}.\frac{(Log\, s)^2}{2} + \xi(b).s^{\lambda-1}.Log\,s +  (\eta_0 + \eta_1.b).s^{\lambda-3} + \zeta(b).s^{\lambda-1} $$
  o\`u \ $\xi, \zeta \in \C[[b]]$ \ et \ $\eta_0, \eta_1 \in \C, \eta_0 \not= 0 $. Alors le th\`eme \ $[\lambda]-$primitif de rang \ $3$ \  \ $\A.e$ \ est isomorphe au quotient 
 $$\A \big/\A.(a - \lambda.b).(1 + u.b + \alpha.b^2)^{-1}.(a - (\lambda+1).b).(a - \lambda.b) $$
 o\`u les nombres complexes \ $u$ \ et \ $\alpha$ \ sont donn\'es par 
 $$ 4\eta_0.u = \eta_1  \quad {\rm et} \quad  4\eta_0.\alpha = (\lambda-1).(\lambda-2).$$
  \end{prop}
 
 \parag{preuve} Calculons d\'ej\`a \ $ X : = (a - \lambda.b).s^{\lambda-1}.\frac{(Log\, s)^2}{2}$.
 \begin{align*}
 & X = s^{\lambda}.\frac{(Log\, s)^2}{2} - \lambda.\big(\frac{s^{\lambda}}{\lambda}.\frac{(Log\, s)^2}{2} - b(\frac{s^{\lambda-1}}{\lambda}.Log\,s)\big) \\
 & X = b(s^{\lambda-1}.Log\,s) = \frac{s^{\lambda}}{\lambda}.Log\,s - \frac{s^{\lambda}}{\lambda}.
  \end{align*}
  Posons \ $Y : = (a - \lambda.b).\xi(b).s^{\lambda-1}.Log\,s $. On a alors
  \begin{align*}
  & Y = \xi(b).(a - \lambda.b).s^{\lambda-1}.Log\,s + b^2.\xi'(b).s^{\lambda-1}.Log\,s \\
  & \qquad = \xi(b).\frac{s^{\lambda}}{\lambda} +  b^2.\xi'(b).s^{\lambda-1}.Log\,s
  \end{align*}
  Calculons \'egalement \ $Z : = (a - \lambda.b).(\eta_0 + \eta_1.b).s^{\lambda-3}$ : 
  \begin{align*}
  & Z = \eta_0.\big[\lambda-2 -\lambda\big].\frac{s^{\lambda-2}}{\lambda-2} + \eta_1.\big[\lambda-1-\lambda\big].\frac{s^{\lambda-1}}{(\lambda-2)(\lambda-1)} \\
 & \qquad = -2\eta_0.\frac{s^{\lambda-2}}{\lambda-2} - \eta_1.\frac{s^{\lambda-1}}{(\lambda-2)(\lambda-1)}
  \end{align*}
  On aura donc
  \begin{align*}
  & (1 + b.\xi'(b))^{-1}.(a - \lambda.b).e = \frac{s^{\lambda}.Log\,s}{\lambda} -\frac{s^{\lambda}}{\lambda} -\frac{2\eta_0}{\lambda-2}.(1 + b.\xi'(b))^{-1}.s^{\lambda-2} + \\
  & \qquad  - \frac{\eta_1}{(\lambda-1)(\lambda-2)}.(1 + b.\xi'(b))^{-1}.s^{\lambda-1} + (1 + b.\xi'(b))^{-1}.(\xi(b) + b.\zeta'(b)).\frac{s^{\lambda}}{\lambda}
  \end{align*}
  On applique alors \ $(a - (\lambda+1).b)$ \ ce qui donne un \'el\'ement de \ $\C[[s]].s^{\lambda-1}$ \ dont  les coefficients respectifs  \ $u,v,w$ \  de \ $s^{\lambda-1}, s^{\lambda}, s^{\lambda+1}$ \ sont respectivement donn\'es par
  \begin{align*}
  & u : =  \frac{4\eta_0}{(\lambda-2)(\lambda-1)} \\
  & v : =  -\frac{4\eta_0.\xi'(0)}{(\lambda-2)(\lambda-1)\lambda} + \frac{\eta_1}{(\lambda-2)(\lambda-1)\lambda} \\
  & w : =  \frac{1}{\lambda+1} - \frac{4\eta_0.(\xi''(0) + 2\xi'(0)^2)}{2(\lambda-2)(\lambda-1)\lambda(\lambda+1)} + \\
 &  \quad - \frac{\eta_1.\xi'(0)}{(\lambda-2)(\lambda-1)\lambda(\lambda+1)} + \frac{2\eta_0.(\xi''(0) + 2\xi'(0)^2)}{(\lambda-2)(\lambda-1)\lambda(\lambda+1)} + \\
 & \quad\quad   \frac{\eta_1.\xi'(0)}{(\lambda-2)(\lambda-1)\lambda(\lambda+1)} =  \frac{1}{\lambda+1}
 \end{align*}
 Donc on trouve
 $$ \frac{4\eta_0}{(\lambda-2)(\lambda-1)}.\big[1 + \frac{\eta_1 - 4\eta_0.\xi'(0)}{4\eta_0}.b + \frac{(\lambda-2)(\lambda-1)}{4\eta_0}.b^2 + \cdots\big].s^{\lambda-1}.$$
 Posons \ $S \in \C[[b]]$ \ l'\'el\'ement entre crochets dans la formule pr\'ec\'edente. On aura donc 
 $$ (a- \lambda.b).S(b)^{-1}.(a - (\lambda+1).b).(1 + b.\xi'(b))^{-1}.(a - \lambda.b).e = 0 .$$
 Posons \ $e_3 : = e, e_2 : = (1 + b.\xi'(b))^{-1}.(a - \lambda.b).e, e_1 : = S(b)^{-1}.(a - (\lambda+1).b).e_2$.\\
 Cherchons maintenant \ $U, V \in \C[[b]]$ \ de mani\`ere que l'\'el\'ement
 $$ \varepsilon_3 : = e_3 + U.e_2 + V.e_1 $$
 v\'erifie
 $$ (a - \lambda.b).\varepsilon_3 = e_2 + (\rho + \sigma.b).e_1.$$
 Cela donne l'\'equation
 \begin{align*}
 & (a - \lambda.b).\varepsilon_3 = e_2 + (b.\xi'(b) + b^2.U'(b) + b.U(b)).e_2 + (U(b).S(b) + b^2.V'(b)) \\
 & =  e_2 + (\rho + \sigma.b).e_1
 \end{align*}
 On veut donc \ $b^2.U'(b) + b.U(b) + b.\xi'(b) = 0$ \ ainsi que \ $\rho = U(0)$ \ et \ $ \sigma = U'(0) + U(0).S'(0)$.\\
 Comme \ $(a - (\lambda+1).b).b.e_1 = 0 $, la valeur du nombre \ $\sigma$ \  ne jouera aucun r\^ole, et on aura \ $\rho = - \xi'(0)$. Alors
 \begin{align*}
 & (a - \lambda).\varepsilon_3 = e_2 + (\rho + \sigma.b).e_1 : = \varepsilon_2 \\
 & (a - (\lambda+1).b).\varepsilon_2 = (S(b) - \rho.b).e_1
 \end{align*}
 On obtient donc que \ $\A.e \simeq \A.\varepsilon_3$ \ est isomorphe \`a 
 $$ \A \big/\A.(a - \lambda.b).(1 + u.b + \alpha.b^2)^{-1}.(a - (\lambda+1).b).(a - \lambda.b) $$
 avec \ $u$ \ et \ $\alpha$ \ les coefficients respectifs de \ $b$ \ et \ $b^2$ \ dans \ $S(b) +\xi'(0).b$, ce qui donne
 $$ u = \frac{\eta_1 - 4\eta_0.\xi'(0)}{4\eta_0} + \xi'(0), \quad \alpha = \frac{(\lambda-2)(\lambda-1)}{4\eta_0} .$$
 On a donc 
  $$u = \frac{\eta_1}{4\eta_0}\qquad {\rm et} \qquad  \alpha = \frac{(\lambda-2)(\lambda-1)}{4\eta_0},$$ ce qui prouve notre assertion.$\hfill \blacksquare$\\

 Prenons maintenant \ $\xi = \zeta = 0$. On a donc
 $$ e : = s^{\lambda-1}.\frac{(Log\, s)^2}{2}  +  (\eta_0 + \eta_1.b).s^{\lambda-3} $$
 et 
 $$\A.e \simeq \A \big/\A.(a - \lambda.b).(1 + u.b + \alpha.b^2)^{-1}.(a - (\lambda+1).b).(a - \lambda.b) $$
 o\`u les nombres complexes \ $u$ \ et \ $\alpha$ \ sont donn\'es par 
 $$ 4\eta_0.u = \eta_1  \quad {\rm et} \quad  4\eta_0.\alpha = (\lambda-1).(\lambda-2).$$
 Posons alors \ $s : = t.(1 + \sigma.t) $ \ o\`u \ $\sigma \in \C$. On obtient
 \begin{align*}
 &  e = t^{\lambda-1}.(1 + \sigma.t)^{\lambda-1}.\frac{[Log\,t + Log(1+\sigma.t)]^2}{2} + \eta_0.t^{\lambda-3}.(1+\sigma.t)^{\lambda-3} + \\
 & \qquad  + \eta_1\frac{t^{\lambda-2}.(1+\sigma.t)^{\lambda-2}}{\lambda-2} 
 \end{align*}
 ce que l'on peut r\'eecrire

 $$ e = S_1(t).t^{\lambda-1}.\frac{(Log\,t)^2}{2} + S_2(t).t^{\lambda}.Log\,t + \theta_0.t^{\lambda-3} + \theta_1.\frac{t^{\lambda-2}}{\lambda-2} + S_3(t).t^{\lambda-1} $$
 o\`u l'on a 
 \begin{align*}
 & S_(t) = 1 + \sigma.(\lambda-1).t + 0(t^2)  \\
 & S_2(t) = 2\sigma + 0(t) \\
 & \theta_0 = \eta_0 \quad {\rm et} \quad  \theta_1 = \eta_1 + \eta_0.(\lambda-2)(\lambda-3).\sigma  \\
 \end{align*}
 la valeur de \ $S_3$ \ importera peu.
 Si on pose \ $\varepsilon : = S_1(t)^{-1}.e $ \ on aura
 \begin{align*}
 & \varepsilon = t^{\lambda-1}.\frac{(Log\, s)^2}{2} + \xi(\beta).t^{\lambda-1}.Log\,t + \eta_0.t^{\lambda-3} +\\
& \quad +  (\theta_1 -\sigma.(\lambda-1)(\lambda-2).\eta_0).\frac{t^{\lambda-2}}{\lambda-2} + \zeta(\beta).t^{\lambda-1}
 \end{align*}
 On a \ $ \theta_1 -\sigma.(\lambda-1)(\lambda-2).\eta_0 = \eta_1 -2\eta_0.\sigma.(\lambda-2) $ \ ce qui montre que les invariants du chang\'e de variables sont \  $4\eta_0\tilde{u} = \eta_1 - 2\eta_0.\sigma.(\lambda-2) $ \ et \ $4\eta_0\tilde{\alpha} = (\lambda-1)(\lambda-2)$ \ et donc \ $\tilde{\alpha} = \alpha$.\\
 On constate donc que pour \ $\sigma \not= 0$ \ le chang\'e de variable n'est pas isomorphe au th\`eme initial, m\^eme si son param\`etre principal lui, n'a pas chang\'e, conform\'ement au th\'eor\`eme \ref{th. fond.}. \\
 On remarquera que le th\`eme de rang 3 consid\'er\'e a la propri\'et\'e d'unicit\'e, c'est-\`a-dire le fait que le couple \ $(\alpha, u)$ \ d\'etermine sa classe d'isomorphie. En effet  \ $E$ \ n'est pas stable car il est de rang 3 et \ $p_2 = 0$ \ et le quotient \ $E\big/F_1$ \  v\'erifie la propri\'et\'e d'unicit\'e\footnote{c'est toujours le cas en rang 2.}  et \ $End_{\A}(E) \simeq \C.\id$ \ puisqu'un endomorphisme de rang 1 devrait envoy\'e \ $e_{\lambda}$ \ sur un \'el\'ement non nul de \ $E_{\lambda+1}$. On peut alors appliquer la proposition 3.3.6 de [B.09b] pour conclure. 

\section{Le cas d'une famille holomorphe.}

Le but de ce paragraphe est de montrer la stabilit\'e de la notion de famille holomorphe de th\`emes \ $[\lambda]-$primitifs par un changement de variable. On obtiendra m\^eme cette stabilit\'e dans le cas relatif, c'est-\`a-dire quand le changement de variable d\'epend holomorphiquement du param\`etre de la famille holomorphe.

Commen{\c c}ons par une formule :

\begin{lemma}\label{formule}
Pour chaque \ $n \in \mathbb{N}$ \ on a
$$ a^n.b = \sum_{p=1}^{n+1}\  \frac{n!}{(n-p+1)!}\, b^p.a^{n-p+1} .$$
\end{lemma}

\parag{Preuve} Par r\'ecurrence sur \ $n\geq 0$. Les cas \ $n=0$ \ et \ $n=1$ \ sont triviaux. Supposons la formule prouv\'ee pour \ $n$. On a alors
\begin{align*}
&  a^{n+1}.b =  \sum_{p=1}^{n+1}\  \frac{n!}{(n-p+1)!} \, a.b^p.a^{n-p+1} \\
& \qquad  =   \sum_{p=1}^{n+1}\  \frac{n!}{(n-p+1)!} \, (b^p.a + p.b^{p+1}).a^{n-p+1} \\
& \qquad  = \sum_{q =1}^{n+2} \ \frac{(n+1)!}{(n+1-q+1)!} \, b^q.a^{n-q+2}
\end{align*}
puisque 
$$ \frac{n!}{(n-q+1)!} + \frac{n!(q-1)}{(n-q+2)!} = \frac{(n+1)!}{(n-q+2)!} . \qquad \qquad \qquad \blacksquare$$

\bigskip

\begin{lemma}
Soit \ $\theta \in \C[[a]]$ \ un changement de variable. On note \ $\alpha : = \theta(a)$ \ et \ $\beta : = b.\theta'(a)$. Pour chaque \ $S \in \C[[b]]$ \ il existe une suite \ $(S_l)_{l\in \mathbb{N}}$ \ de \ $\C[[b]]$ \ telle que l'on ait
$$ S(b) = \sum_{\l \geq 0} \ S_l(\beta).\alpha^l $$
dans l'alg\`ebre \ $\hat{A}$.
\end{lemma}

\parag{Preuve} Il suffit de montrer que pour chaque \ $\nu \in \mathbb{N}$ \ on peut trouver une suite \ $(S_{\nu,l})_{l \geq 0}$ \ dans \ $\C[[\beta]]$ \ de mani\`ere \`a v\'erifier
$$ b^{\nu} = \beta^{\nu}.\sum_{\l \geq 0} \ S_{\nu,l}(\beta).\alpha^l  \quad {\rm avec} \quad S_{\nu,0}(0) = \theta'(0)^{-\nu}.$$
Ceci s'obtient facilement par r\'ecurrence sur \ $\nu$. Posons \ $a = \eta(\alpha)$, c'est \`a dire que \ $\eta : = \theta^{-1}$ \ au sens de la composition des s\'eries formelles sans terme constant. Si l'on a
$$ b^{\nu} = \beta^{\nu}.\sum_{\l \geq 0} \ S_{\nu,l}(\beta).\alpha^l $$
on aura
\begin{equation*}
 b^{\nu+1} =  \beta^{\nu}.\big[\sum_{\l \geq 0} \ S_{\nu,l}(\beta).\alpha^l \big].\beta.\eta'(\alpha)
\end{equation*}
et il suffit d'appliquer le lemme \ref{formule} pour faire avancer la r\'ecurrence, en constatant que \ $\eta'(0) = \theta'(0)^{-1}$. $\hfill \blacksquare$
 
 \begin{prop}\label{chgt thematique}
 Soit \ $X$ \ un espace complexe r\'eduit et soit \ $\varphi : X \to \Xi_{\lambda}^{(N-1)}$ \ une application \ $k-$th\'ematique. Soit \ $\theta \in \C[[a]]$ \ un changement de variable. Soit \ $\theta_*(\mathbb{E}_{\varphi})$ \ le faisceau de \ $\mathcal{O}_X-$modules \ $\mathbb{E}_{\varphi}$ \ muni des op\'erations \ $\alpha : = \theta(a)$ \ et \ $\beta : = b.\theta'(a)$. Alors \ $\theta_*(\mathbb{E}_{\varphi})$ \ est un \ $\mathcal{O}_X[[\beta]]-$module libre de rang \ $k$ \ de base \ $\varphi, \alpha.\varphi, \dots, \alpha^{k-1}.\varphi$ \ et qui est stable par \ $\alpha$.\\
 Donc \ $\varphi$ \ est encore \ $k-$th\'ematique pour la \ $\mathcal{O}_X(\alpha,\beta)-$structure de \ $\Xi_{\lambda}^{(N-1)}$ \ d\'efinie par \ $\alpha$ \ et \ $\beta$.
 \end{prop}
 
 \parag{Preuve} Il s'agit de montrer que \ $\varphi, \alpha.\varphi, \dots, \alpha^{k-1}.\varphi$ \ est une \ $\mathcal{O}_X[[\beta]]-$base de \ $\mathcal{E}_{\varphi}$, ce faisceau \'etant manifestement stable par \ $\alpha$.\\
 Comme on a \ $b.\mathbb{E} = \beta.\mathbb{E}$, et \ $a^k.\mathbb{E} \subset b.\mathbb{E}$, la matrice de \ $\varphi, \alpha.\varphi, \dots, \alpha^{k-1}.\varphi$ \ dans la \ $\mathcal{O}_X-$base \ $\varphi, a.\varphi, \dots, a^{k-1}.\varphi$ \ de \ $\mathbb{E}\big/b.\mathbb{E}$ \ est triangulaire avec des \ $1$ \ sur la diagonale. Ceci montre que \ $\varphi, \alpha.\varphi, \dots, \alpha^{k-1}.\varphi$ \ est une une \ $\mathcal{O}_X-$base de  \ $\mathbb{E}\big/b.\mathbb{E} \simeq \mathbb{E}\big/\beta.\mathbb{E}$. On a donc bien une \ $ \mathcal{O}_X[[\beta]]-$base. Comme, de plus, on a \ $\alpha^k.\varphi \in b.\mathbb{E} = \beta.\mathbb{E}$, on a bien montr\'e que \ $\varphi$ \ est \ $k-$th\'ematique pour la (a,b)-structure donn\'ee par \ $(\alpha, \beta)$. $\hfill \blacksquare$\\
 
 \begin{defn}\label{chgt avec param.}
 Soit \ $X$ \ un espace complexe r\'eduit. Nous dirons que l'\'el\'ement  \ $\theta \in \mathcal{O}(X)[[a]]$ \ est un {\bf changement de variable \ $X-$relatif} \ si l'on a \ $\theta(0) = 0$ \ et \ $\theta'(0)(x) \not= 0$ \ pour tout \ $x \in X$.
 \end{defn}
 
 On a donc pour chaque \ $x \in X$ \ fix\'e un changement de variable.
 
 Une cons\'equence imm\'ediate de la  proposition \ref{chgt thematique} est le th\'eor\`eme suivant de stabilit\'e par changement de variable relatif  pour les familles holomorphes de th\`emes \ $[\lambda]-$primitifs de rang \ $k$.\\
 
 \begin{thm}\label{chgt famille}
 Soit \ $X$ \ un espace complexe r\'eduit et soit \ $\mathbb{E}$ \ une famille holomorphe de th\`emes \ $[\lambda]-$primitifs de rang \ $k$ \ param\'etr\'ee par \ $X$. Consid\'erons un changement de variable  \ $X-$relatif \ $\theta \in \mathcal{O}(X)[[a]]$.  Alors la famille \ $\theta_*(\mathbb{E})$ \ des \ $\theta_*(\mathbb{E}(x))_{x \in X}$, est holomorphe.
 \end{thm}
 
 Donnons  un  corollaire \'egalement  imm\'ediat, mais important de ce th\'eor\`eme et du th\'eor\`eme \ref{th. fond.}.
 
 \begin{cor}\label{inv. param. principaux}
 Soit \ $X$ \ un espace complexe r\'eduit connexe et soit \ $\mathbb{E}$ \ une famille holomorphe de th\`emes \ $[\lambda]-$primitifs de rang \ $k$ \ param\'etr\'ee par \ $X$. Soit \ $\theta \in \mathcal{O}(X)[[a]]$ \ un changement de variable relatif. Notons \ $(\lambda_1, p_1, \dots, p_{k-1})$ \ les invariants fondamentaux de cette famille, et notons \ $\alpha_1, \dots, \alpha_{k-1} : X \to \C$ \ les fonctions holomorphes sur \ $X$ \ donn\'ees par les param\`etres principaux de cette famille.\\
 Alors les fonctions holomorphes sur \ $X$ \ donn\'ees par les param\`etres principaux de la famille holomorphe \ $\theta_*(\mathbb{E})$ \ sont respectivement \ $r^{p_1}.\alpha_1, \dots, r^{p_{k-1}}.\alpha_{k-1}$, o\`u l'on a pos\'e  \ $r : = \theta'(0)\in \mathcal{O}(X)$.
 \end{cor}

 \bigskip
 
 \bigskip

\section{R\'ef\'erences}

\begin{itemize}

\item {[B.93]} Daniel Barlet {\it Th\'eorie des (a,b)-modules I}, Univ. Ser. Math. Plenum (1993), p.1-43.

\item {[B.08]} Daniel Barlet {\it Sur certaines singularit\'es d'hypersurfaces II}, J. Alg. Geom. 17 (2008), p. 199-254.

\item {[B.09a]} Daniel Barlet {\it P\'eriodes \'evanescentes et (a,b)-modules monog\`enes}, Bollettino U.M.I (9) II (2009), p.651-697.

\item {[B.09b]} Daniel Barlet {\it Le th\`eme d'une p\'eriode \'evanescente}, preprint de l'Institut E. Cartan (Nancy) 2009  $n^0$  33,  57 pages.

\end{itemize}

\end{document}